\documentclass[11pt,fleqn]{article}
\usepackage{latexsym,amsfonts,amsbsy,graphics,epsf}
\newtheorem{lemma}{Lemma}
\newtheorem{corollary}{Corollary}
\newtheorem{theorem}{Theorem}
\newtheorem{definition}{Definition}

\def\Integer{\mathbb{Z}}
\def\Mult#1#2{\left[{#1 \atop #2}\right]}
\def\Mults#1#2{\left[{#1 \atop #2}\right]'}
\def\mod#1{\!\!\! \pmod{#1}}
\def\la{\lambda}
\def\vi{\boldsymbol i}
\def\vv{\boldsymbol v}
\def\vet{\boldsymbol \eta}
\def\vn{\boldsymbol n}
\def\vm{\boldsymbol m}
\def\vmu{\boldsymbol \mu}
\def\ve{\boldsymbol e}
\def\bp{$(\alpha,\beta)$}
\def\bpp{$(\alpha',\beta')$}
\def\gd{$(\gamma,\delta)$}
\def\GD{$(\Gamma,\Delta)$}

\def\proof#1{{\em Proof.} #1 \hfill $\Box$}
\def\proofarg#1#2{{\em #1} #2 \hfill $\Box$}
\def\sumv#1{{\sum_{#1 \in \Integer^{N-1}}}^{\hspace{-3mm}(\sigma)}}
\def\sumvb#1{{\sum_{#1 \in \Integer^{N-1}}}^{\hspace{-4mm}(\sigma)}}
\def\sv#1#2{{\sum_{#1 \in \Integer^{N-1}}}^{\hspace{-3mm}(\sigma,#2)}}
\def\svb#1#2{{\sum_{#1 \in \Integer^{N-1}}}^{\hspace{-4mm}(\sigma,#2)}}

\def\nsection#1{\setcounter{equation}{0}\setcounter{lemma}{0}
\setcounter{theorem}{0}\setcounter{corollary}{0}\setcounter{definition}{0}
\section{#1}}

\begin{document}


\title{A Higher-Level Bailey Lemma: \\
Proof and Application}

\author{\Large
Anne Schilling\thanks{
e-mail: {\tt anne@insti.physics.sunysb.edu}} \\
\mbox{} \\
\em Institute for Theoretical Physics, State University of New York \\
\em Stony Brook, NY 11794-3840, USA \\
\mbox{} \\
\Large and \\
\mbox{} \\
\Large S.~Ole Warnaar\thanks{
e-mail: {\tt warnaar@maths.mu.oz.au}} \\
\mbox{} \\
\em Mathematics Department, University of Melbourne\\
\em Parkville, Victoria 3052, Australia}

\date{\Large July, 1996}
\maketitle
\begin{abstract}
In a recent letter, new representations were proposed for the pair of 
sequences ($\gamma,\delta$), as defined formally by Bailey in his 
famous lemma. Here we extend and prove this result, providing pairs 
($\gamma,\delta$) labelled by the Lie algebra A$_{N-1}$, two non-negative 
integers $\ell$ and $k$ and a partition $\lambda$, whose parts do not 
exceed $N-1$. Our results give rise to what we call a ``higher-level'' 
Bailey lemma. As an application it is shown how this lemma can be applied 
to yield general $q$-series identities, which generalize some 
well-known results of Andrews and Bressoud.
\end{abstract}

\nsection{Introduction}
A well-known approach to proving $q$-series identities of the
Rogers--Ramanujan type is the Andrews--Bailey 
construction~\cite{Bailey,Andrews84}, given by the following lemma.
\begin{lemma}\label{LAB}
Let $\alpha=\{\alpha_L\}_{L\geq 0}$ and $\beta=\{\beta_L\}_{L\geq 0}$
be sequences which satisfy
\begin{equation}
\beta_L = \sum_{r=0}^L \frac{\alpha_r}{(q)_{L-r}(aq)_{L+r}} \, ,
\label{ab}
\end{equation}
for all $L$.
Then the sequences $\alpha'$ and $\beta'$ defined by
\begin{eqnarray}
\alpha'_L &=&
\frac{(\rho_1)_L (\rho_2)_L (aq/\rho_1 \rho_2)^L}
{(aq/\rho_1)_L (aq/\rho_2)_L} \: \alpha_L \nonumber \\[2mm]
\beta'_L &=& \sum_{r=0}^L
\frac{(\rho_1)_r (\rho_2)_r (aq/\rho_1 \rho_2)^r
(aq/\rho_1 \rho_2)_{L-r}}
{(aq/\rho_1)_L (aq/\rho_2)_L (q)_{L-r}} \: \beta_r ,
\label{AB}
\end{eqnarray}
again satisfy equation (\ref{ab}).
\end{lemma}
A pair of sequences \bp \ that satisfies (\ref{ab}) is called
a Bailey pair relative to $a$. The iteration of lemma~\ref{LAB}:
$$(\alpha,\beta) \to (\alpha',\beta') \to (\alpha'',\beta'') \to \cdots $$
is called a Bailey chain.
The symbol $(a)_n$ in the above lemma is the usual $q$-shifted
factorial,
\begin{equation}
(a;q)_{\infty}=(a)_{\infty} = \prod_{n=0}^{\infty} (1-aq^n)
\end{equation}
and
\begin{equation}
(a;q)_n = (a)_n = \frac{(a;q)_{\infty}}{(aq^n;q)_{\infty}}
\qquad {n \in \Integer}.
\end{equation}

Perhaps not so well-known however, is that Bailey in his original
paper~\cite{Bailey} considered another pair of sequences
$\gamma=\{\gamma_L\}_{L\geq 0}$ and $\delta=\{\delta_L\}_{L\geq 0}$ 
related by
\begin{equation}
\gamma_L = \sum_{r=L}^\infty \frac{\delta_r}{(q)_{r-L}(aq)_{r+L}} \, ,
\label{gd}
\end{equation}
for all $L$.
We will refer to a pair of sequences \gd \ which
satisfies (\ref{gd}) as a conjugate Bailey pair relative to $a$.
Bailey noted that given a Bailey pair and some conjugate Bailey pair
(both relative to $a$) the following identity holds:
\begin{equation}
\sum_{L=0}^{\infty} \alpha_L \gamma_L =
\sum_{L=0}^{\infty} \beta_L \delta_L \, .
\label{agbd}
\end{equation}
As pointed out by Andrews~\cite{Andrews84}, 
lemma~\ref{LAB} is a special case of this 
result with conjugate Bailey pair
\begin{equation}
\gamma_L = \frac{(\rho_1)_L (\rho_2)_L (aq/\rho_1\rho_2)^L}
{(aq/\rho_1)_L (aq/\rho_2)_L} \:
\frac{1}{(q)_{M-L}(aq)_{M+L}}
\label{gammarho}
\end{equation}
and
\begin{equation}
\delta_L = \frac{(\rho_1)_L (\rho_2)_L (aq/\rho_1\rho_2)^L}
{(aq/\rho_1)_M (aq/\rho_2)_M} \:
\frac{(aq/\rho_1\rho_2)_{M-L}}{(q)_{M-L}} \, ,
\label{deltarho}
\end{equation}
where $M$ is an arbitrary non-negative integer.

As a special case of (\ref{gammarho}) and (\ref{deltarho})
we may send the undetermined parameters $\rho_1$ and $\rho_2$ to infinity.
Using 
\begin{equation}
\lim_{a\to \infty} a^{-n} (a)_n = (-1)^n q^{n(n-1)/2} ,
\end{equation}
this gives
the simple expressions
\begin{equation}
\gamma_L = \frac{a^L q^{L^2}}{(q)_{M-L}(aq)_{M+L}}
\qquad {\rm  and } \qquad
\delta_L = \frac{a^L q^{L^2}}{(q)_{M-L}} \, .
\label{gammadelta}
\end{equation}
It is the aim of this paper to show that these two expressions
are the first instances of an infinite series of conjugate Bailey
pairs.
Specifically, we establish a conjugate Bailey pair relative to 
$a=q^{k+\ell}$, for each pair of non-negative integers $k$ 
and $\ell$,
positive integer $N$ and partition $\la$ with largest
part not exceeding $N-1$.
For $k=0$ and $\la=(\ell)$ or $|\la|=0,~\ell=N$, this result 
was announced in ref.~\cite{SW}.
The conjugate Bailey pair (\ref{gammadelta}) with $a$ of the form
$a=q^{\ell}$ is the special case $N=1$ and $k=0$.

Before we proceed we note that Milne and Lilly~\cite{LM93,ML92,ML95} 
have recently given higher-dimensional generalizations of the Bailey 
transform and the Bailey lemma in the setting of very-well poised 
basic hypergeometric functions on the higher rank groups A$_{\ell}$ 
and C$_{\ell}$. 
Specifically they generalized the defining relation
of a Bailey pair to higher rank (equation  (\ref{LAB}) corresponding
to A$_1$) and obtained the corresponding generalization
of lemma~\ref{LAB}.
This has to be contrasted with the results obtained
in this paper, which stay within the framework of A$_1$,
but provide higher-level (spin) representations for conjugate Bailey
pairs.
Indeed these representations involve expressions closely related
to the string functions of the affine Lie algebra A$_1^{(1)}$ at level
$N$, and can be used to obtain 
$q$-series identities for the branching functions of the level-$N$ 
A$_1^{(1)}$ cosets theories.  

In the next section we slightly reformulate the definition of a 
conjugate Bailey pair and state our main result; an infinite hierarchy
of conjugate Bailey pairs.
In section~\ref{secproof} we present a recursive proof of our result using
telescopic expansion techniques. In section~\ref{secapp} we apply what
we term the ``higher-level'' Bailey lemma to derive a very
general $q$-series identity, which contains many 
known Rogers-Ramanujan type identities as specializations. 
For the proof of this general identity we also
introduce some further transformations of the Bailey chain and
Bailey lattice type.
We finally conclude with a brief discussion of our results.

\nsection{A hierarchy of conjugate Bailey pairs}\label{secCBP}
In the following we present an infinite series of pairs of
sequences \gd \ which satisfy equation (\ref{gd}).
In fact, for notational reasons it turns out to be 
fruitful to slightly modify the definition (\ref{gd})
of a conjugate Bailey pair.
\begin{definition}
A pair of sequences 
$\Gamma=\{\Gamma_{L,k}\}_{L,k\geq 0}$ and 
$\Delta=\{\Delta_{L,k}\}_{L\geq k\geq 0}$
which satisfies
\begin{equation}
\Gamma_{L,k} = \sum_{r=L+k}^\infty 
\frac{\Delta_{r,k}}{(q)_{r-L-k}(aq)_{r+L}} \, ,
\label{CD}
\end{equation}
for all $L,k\geq 0$ is called a \GD-pair relative to $a$.
\end{definition}
At first sight it may seem that a \GD-pair is more general
than a conjugate Bailey pair, reducing to it for $k=0$.
However note that given a \GD-pair relative to $a$ and setting 
\begin{equation}
\gamma_L^{(k)} = \Gamma_{L,k} 
\qquad {\rm  and } \qquad
\delta_L^{(k)} = \frac{\Delta_{L+k,k}}{(aq)_k} \, ,
\label{CDcd}
\end{equation}
the pair $(\gamma^{(k)},\delta^{(k)})$ becomes a 
conjugate Bailey pair relative to $aq^k$.

\vspace{3mm}
We now proceed to state a series of \GD-pairs.
For this some further notation is needed.
First, we need two types of Gaussian polynomials or $q$-binomial
coefficients
(see for example \cite{GR}).
\begin{definition}
For $n,m$ arbitrary integers 
\begin{equation}
\Mult{m+n}{n} = \left\{
\begin{array}{ll}
\displaystyle \frac{(q)_{n+m}}{(q)_n (q)_m} \quad & n,m \geq 0 \\[4mm]
0 & {\rm otherwise,}
\end{array} \right.
\label{Gpoly}
\end{equation}
and
\begin{equation}
\Mults{m+n}{n}=
\left\{ \begin{array}{ll}
\displaystyle
\frac{(q^{n+1})_m}{(q)_m} \quad & m\geq 0 \\[4mm]
0 & {\rm otherwise.} \end{array} \right.
\label{bins}
\end{equation}
\end{definition}
Notice that the primed version of the $q$-binomials is asymmetric in 
$m$ and $n$, and that both types of $q$-binomials coincide for $m+n\geq 0$.
Both the ordinary and the primed binomials enjoy the recurrences
\begin{eqnarray}
\Mult{m+n}{n}&=&\Mult{m+n-1}{n-1}+q^n\Mult{m+n-1}{n}
\label{arr}
\\[2mm]
&=&\Mult{m+n-1}{n}+q^m\Mult{m+n-1}{n-1},
\label{brr}
\end{eqnarray}
with the notable exception that when $n=m=0$, (\ref{arr}) and (\ref{brr})
break down for the ordinary $q$-binomials (\ref{Gpoly}).

Next we introduce the following notation.
Throughout $N$ is assumed to be a fixed positive integer
and $\vv=(v_1,\ldots,v_{N-1})$ denotes
a vector with integer entries, i.e., $\vv\in \Integer^{N-1}$.
In particular, $\ve_j$ is the $j$-th unit vector,
$(\ve_j)_k=\delta_{j,k}$ ($\delta_{j,k}=1$ for $j=k$ and 0 otherwise).
By definition $\ve_0=\ve_{N}=\boldsymbol{0}$.
For a given matrix $M$, we use the following notation for matrix 
multiplication $\vv M \vv = \sum_{j,k=1}^{N-1} v_j M_{j,k} v_k $, 
where we omit the transposition symbol. 
We will often encounter the Cartan matrix of the Lie
algebra A$_{N-1}$, denoted by $C$. 
That is, $C_{j,k}=2\delta_{j,k}-\delta_{j,k-1}-\delta_{j,k+1}$.
Further we adopt the compact notation
\begin{equation}
\Mult{k}{\vv}_q =\frac{(q)_k}{(q)_{v_1} \cdots (q)_{v_{N-1}}
(q)_{k-v_1-\cdots-v_{N-1}}}
\end{equation}
for $\vv\in \Integer^{N-1}_{\geq 0}$ and zero otherwise.
Lastly, let $\la=(\la_1,\la_2,\ldots)$, $\la_1\geq \la_2 \geq \ldots$ denote
a partition with $\sum_i \la_i = |\la|$.
If $\la$ has largest part $\la_1\leq N-1$, 
we abbreviate $\ve_{\la_1}+\ve_{\la_2}+\ldots$ to $\ve_{\la}$.

We can now state the following theorem.
\begin{theorem}\label{thGD}
Fix integers $M\geq 0$, $N\geq 1$ and $\ell\geq 0$, and fix 
a partition $\la$ whose largest part does not exceed $N-1$.
Choose $\sigma \in \{0,1\}$ such that $\ell+|\la|+\sigma N$ is even.
Then the following two sequences form a
\GD-pair relative to $a=q^{\ell}$:
\begin{eqnarray}
\lefteqn{
\Gamma_{L,k}=\frac{a^{L/N+k} q^{L^2/N+kL}}{(q)_{M-L-k}(aq)_{L+M}}}
\nonumber \\[2mm]
& & \hspace{1cm} \times \sumvb{\vet,\vi}
q^{-\vi\, C^{-1}\bigl(\vi+(2L+\ell)\ve_1\bigr)}
\Mult{k}{\vi}_{\! 1/q} 
q^{\vet \, C^{-1}(\vet-\ve_{\la})}\prod_{j=1}^{N-1}
\Mults{\mu_j+\eta_j}{\eta_j}
\label{GammaN}
\end{eqnarray}
and
\begin{equation}
\Delta_{L,k} =
\frac{a^{L/N}q^{L^2/N-kL}}{(q)_{M-L}}
\sumv{\vn} q^{\vn\, C^{-1} (\vn-\ve_{\la})}
\prod_{j=1}^{N-1} \Mult{m_j+n_j}{n_j}.
\label{DeltaN}
\end{equation}
The variables $\mu_j$, which occur in the primed $q$-binomial
in (\ref{GammaN}), are fixed by the (external) variables
$L,k$ and the (summation) variables $\vi$ and $\vet$
by the equation
\begin{equation}
\vmu=C^{-1}\Bigl((M-L-k)\ve_1+(M+L+\ell)\ve_{N-1}+\ve_{\la}-
\sum_{j=1}^{N-1} (i_j+i_{j+1})\ve_j-2\vet\Bigr),
\label{muetai}
\end{equation}
where 
\begin{equation}
i_N := k-i_1-\cdots - i_{N-1}.
\label{iN}
\end{equation}
 Similarly, the variables $m_j$, which occur in the $q$-binomial
 in (\ref{DeltaN}), are fixed by $L,\vn$ thanks to
\begin{equation}
\vm = C^{-1}\Bigl((2L+\ell)\,\ve_{N-1}+\ve_{\la}-2\vn\Bigr).
\label{mn}
\end{equation}
\end{theorem}
So far we refrained from defining the
superscript $(\sigma)$ in the sums over $\vet,~\vi$ and $\vn$ in 
(\ref{GammaN}) and (\ref{DeltaN}).
To explain this notation (as well as the origin of the parity restrictions
in the theorem) let us consider expression (\ref{DeltaN}).
In this expression we sum over the vector $\vn$ with integer entries.
Given such $\vn$, we can compute a companion vector $\vm$ using the
$(\vm,\vn)$-system (\ref{mn}). Of course, because of the
product over the $q$-binomials, we only wish to consider
those $\vn$ which return a vector $\vm$ whose entries are all integer.
A close scrutiny of equation (\ref{mn}) shows that 
$\vm\in \Integer^{N-1}$ if and only if 
$\ell+|\la|+\sigma N$ is even and 
$$\frac{L+(\ell-|\la|)/2}{N}-(C^{-1} \vn)_1 \in \Integer+\frac{\sigma}{2},$$
where $\sigma$ takes either the value 0 or 1.
A similar analysis can be made for the sum in (\ref{GammaN}) where $\vn$
is replaced by $\vet+\vi$.
More compactly we denote
\begin{equation}
\sum_{\vn_1,\ldots,\vn_p \in \Integer^{N-1} \atop
\frac{L+(\ell-|\la|)/2}{N}-(C^{-1} (\vn_1+\cdots+\vn_p))_1 
\in \Integer +\frac{\sigma}{2}}
={\sum_{\vn_1,\ldots,\vn_p \in \Integer^{N-1}}}^{\hspace{-8mm}(\sigma,L)}
={\sum_{\vn_1,\ldots,\vn_p \in \Integer^{N-1}}}^{\hspace{-8mm}(\sigma)}.
\label{sumsymbol}
\end{equation}

\vspace{5mm}
After stating theorem~\ref{thGD}, it seems appropriate
to make a few remarks about its origin.
In statistical mechanics there has recently been much interest
in so-called fermionic representations of Virasoro characters,
see e.g., ref.~\cite{SB}. Such fermionic character
representations relate to combinatorial or ``subtractionless'' 
bases for Virasoro modules. One of the most successful approaches
to this problem has been the study of exactly solvable lattice 
models using what is termed the Bethe Ansatz technique.
One of the many remarkable features of this technique, is that
it gives rise to so-called fermionic polynomials, which in 
a special limit give fermionic character expressions.
(The polynomials featuring in items 287 and 289 of 
MacMahon's Combinatory Analysis, Vol.~2~\cite{MacMahon},
are probably the simplest examples of fermionic polynomials.)
The rich combinatorics of the fermionic polynomials are described
by systems of (quasi)particles that obey exclusion statistics.
These exclusion rules are encoded in what are called 
$(\vm,\vn)$-systems~\cite{Berkovich}.
Here $\vn$ is a vector whose $j$-th entry $n_j$ denotes the occupation
number of particles of type $j$. The vector $\vm$ has an interpretation
in terms of occupation numbers of holes or anti-particles.

Returning to the above theorem, equations (\ref{muetai}) and (\ref{mn})
are examples of the afore-mentioned $(\vm,\vn)$-systems\footnote{For 
the ``interpretation'' of the $(\vmu,\vet)$-system (\ref{muetai}),
one in fact has to allow for the occupation numbers $\eta_j$ to
become negative. A phenomenon discussed in~\cite{BMO,BMS96}.},
and apart from 
an overall factor, the sequence $\Delta$ with $\sigma=0$ and
$\la=(\ell)$, coincides with the
fermionic polynomials of the level-2 A$_{N-1}^{(1)}$ Jimbo--Miwa--Okado
models~\cite{JMO}, as obtained in ref.~\cite{FOW}. 

\vspace{5mm}
Maybe more important than theorem~\ref{thGD} itself is the following 
corollary.
\begin{corollary}\label{corgd}
Let $M,N,\ell,\la$ and $\sigma$ be as in theorem~\ref{thGD}.
Then the following two sequences form a
conjugate Bailey pair relative to $a=q^{\ell}$:
\begin{equation}
\gamma_L =
\frac{a^{L/N}q^{L^2/N}}{(q)_{M-L}(aq)_{M+L}}
\sumv{\vet} q^{\vet \, C^{-1} (\vet-\ve_{\la})}
\prod_{j=1}^{N-1} \Mult{\mu_j+\eta_j}{\eta_j}
\label{gammaN}
\end{equation}
and
\begin{equation}
\delta_L =
\frac{a^{L/N}q^{L^2/N}}{(q)_{M-L}}
\sumv{\vn} q^{\vn\, C^{-1} (\vn-\ve_{\la})}
\prod_{j=1}^{N-1} \Mult{m_j+n_j}{n_j}.
\label{deltaN}
\end{equation}
Here $\vm$ is given by (\ref{mn}) and
$\vmu$ by
\begin{equation}
\vmu=C^{-1}\Bigl((M-L)\ve_1+(M+L+\ell)\ve_{N-1}+\ve_{\la}-2\vet\Bigr).
\label{mueta}
\end{equation}
\end{corollary}
We remark that this generalizes the conjugate
Bailey pair (\ref{gammadelta}) (with the restriction $a=q^{\ell}$,
$\ell \in \Integer_{\geq 0}$.)

\proofarg{Proof of corollary~\ref{corgd}}{\newline
Setting $k=0$ in equation~(\ref{DeltaN}) immediately gives 
$\Delta_{L,0}=\delta_L$.
Setting $k=0$ in equation~(\ref{GammaN}) we find that the
sum over $\vi$ is non-zero for $\vi=\boldsymbol{0}$ only.
Hence (\ref{muetai}) simplifies to (\ref{mueta}) and
$\Gamma_{L,0}$ gives precisely the above expression
for $\gamma_L$ provided one can show that 
the primed $q$-binomials in (\ref{GammaN}) can
be replaced by the unprimed $q$-binomials.
To see this, notice that both $q$-binomials are zero unless
$\mu_j\geq 0$ for all $j=1,\ldots,N-1$. From this, and (\ref{mueta}),
it follows that with $\mu_0=\mu_N=0$,
\begin{equation}
\mu_j+\eta_j=\frac{1}{2}\Bigl(\mu_{j-1}+\mu_{j+1}+(M-L)\delta_{1,j}
+(M+L+\ell)\delta_{N-1,j}+\delta_{\la_1,j}+\delta_{\la_2,j}
+\cdots \Bigr) \geq 0.
\end{equation}
(We can safely assume that $M\geq L$ since $\Gamma_{L,0}$ (and $\gamma_L$)
are non-zero for $M\geq L$ only.)
But for $\mu_j+\eta_j \geq 0$ the primed and unprimed $q$-binomials
coincide, and so indeed $\Gamma_{L,0}= \gamma_L$. }

To conclude this section we note that corollary~\ref{corgd} will
be used in the following for applications rather than theorem~\ref{thGD}.
Almost all of the identities arising from (\ref{GammaN})
and (\ref{DeltaN}) can either be obtained from corollary~\ref{corgd}
or are of such complexity that they seem of little relevance.
Still, the parameter $k$ in theorem~\ref{thGD} is indispensable
to us, as it gives rise to recurrences which are at the heart
of our proof of the infinite series of conjugate Bailey pairs.

\nsection{Proof of theorem~\ref{thGD}}\label{secproof}
In this section we prove theorem~\ref{thGD} using recurrence relations.
For the proof of these recurrences we rely on the following
generalizations of the $q$-binomial recurrences
(\ref{brr}), termed telescopic expansions~\cite{Berkovich}.
(Similar generalizations of (\ref{arr}) will not be needed here.)
\begin{lemma}
For integer $N\geq 2$ and $A_j,B_j$ integer for all $j=1,\ldots,N-1$,
\begin{equation}
\prod_{j=1}^{N-1}\Mults{A_j+B_j}{A_j}
=\prod_{j=1}^{N-1}\Mults{A_j+B_j-1}{A_j}
+\sum_{p=1}^{N-1} q^{B_p}
\prod_{j=1}^{N-1}\Mults{A_j+B_j-\chi(j\leq p)}{A_j-\delta_{j,p}}
\label{rtele}
\end{equation}
and
\begin{equation}
\prod_{j=1}^{N-1}\Mults{A_j+B_j}{A_j}
=\prod_{j=1}^{N-1}\Mults{A_j+B_j-1}{A_j}
+\sum_{p=1}^{N-1} q^{B_p} \prod_{j=1}^{N-1}
\Mults{A_j+B_j-\chi(j\geq p)}{A_j-\delta_{j,p}}
\label{btele}
\end{equation}
where $\chi({\rm true})=1$ and $\chi({\rm false})=0$.
\end{lemma}
\proof{For $N=2$ (\ref{rtele}) and (\ref{btele}) simplify
to (\ref{brr}). For $N\geq 2$, (\ref{rtele}) and (\ref{btele}) follow
by induction on $N$ and application of (\ref{brr}).} 

Let us now return to theorem~\ref{thGD}.
Inserting both the definition of $\Delta$ and $\Gamma$
into equation (\ref{CD})
and reshuffling some of the terms,
the claim of theorem~\ref{thGD} can be
re-expressed as the following polynomial identity.
\begin{lemma}\label{lemGD}
Let $M,N,\ell,\la$ and $\sigma$ be as in theorem~\ref{thGD}.
Then for $0 \leq L \leq M$ and $0\leq k\leq M-L$,
\begin{eqnarray}
\lefteqn{
\svb{\vet,\vi}{L} q^{-\vi\, C^{-1}\bigl(\vi+(2L+\ell)\ve_1\bigr)}
\Mult{k}{\vi}_{\! 1/q}
q^{\vet \, C^{-1}(\vet-\ve_{\la})}\prod_{j=1}^{N-1}
\Mults{\eta_j+\mu_j}{\eta_j} } \nonumber\\[2mm]
&=&
\sum_{r=L+k}^M q^{(r+L+\ell)(r-L-Nk)/N}
\Mult{M-L-k}{M-r}
\frac{(q^{\ell+1})_{L+M}}{(q^{\ell+1})_{L+r}} \nonumber \\[2mm]
& & \hspace{3cm} \times
\sv{\vn}{r} q^{\vn \, C^{-1}(\vn-\ve_{\la})}\prod_{j=1}^{N-1}
\Mult{m_j+n_j}{n_j},
\label{lem}
\end{eqnarray}
with $(\vmu,\vet)$-system (\ref{muetai}) and
$(\vm,\vn)$-system (\ref{mn}), where $L$ in (\ref{mn}) 
has been replaced by $r$.
\end{lemma}

\subsection{Proof of lemma~\ref{lemGD}}
We prove the polynomial identity
(\ref{lem}) by induction on $k$ and $M-L$.
To this end we show that both the left- and right-hand side of (\ref{lem}),
denoted by $f_1$ and $f_2$, respectively,
obey the recursion relations
\begin{eqnarray}
\lefteqn{
f(M,L,k)=f(M-1,L,k) +
q^{M+L+\ell} \Bigl(f(M,L,k+1)-f(M-1,L,k)\Bigr)}  \nonumber \\
& & \hspace{8.43cm} {\rm for} \;  0\leq k<M-L \label{rec} \\
\lefteqn{f(M,L,k)=q^{-(2L+\ell+1)(N-1)/N}f(M,L+1,k-1)
\qquad {\rm for} \;  k=M-L} \label{boundrec}
\end{eqnarray}
as well as the initial condition
\begin{equation}
f_1(M,M,0)=f_2(M,M,0) \qquad {\rm for} \; M\geq 0.
\label{initial}
\end{equation}

Equations (\ref{rec})--(\ref{initial}) fix the function $f$ uniquely.
With (\ref{initial}) as a seed, recurrence 
(\ref{boundrec}) gives $f(M,L,k)$, for all $k=M-L$.
Now assume 
that $f(M,L,k)$ is known for all $k=M-L-n$, with
$n=0,\ldots,n_0$. Since (\ref{rec}) expresses each $f(M,L,k)$ with
$k=M-L-n-1$ as the sum of two $f(M,L,k)$ with $k=M-L-n$,
this implies that $f(M,L,k)$ is fixed for
$k=K-L-(n_0+1)$. Since the assumption is obviously true for $n_0=0$,
we indeed find that $f$ is fixed by
(\ref{rec})--(\ref{initial}).

\begin{figure}
\epsfxsize = 5cm
\centerline{\epsffile{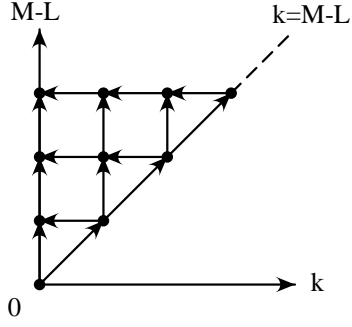}}
\caption{Graphical representation of the ``flow'' induced by the 
recurrence relations (\ref{rec})--(\ref{initial}) 
in the $(k,M-L)$-plane.}\label{fig}
\end{figure}

The recurrences are illustrated in figure \ref{fig}. Each coordinate
in the $(k,M-L)$-plane for which $0\leq L \leq M$ and $0\leq k \leq
M-L$ is represented by a node,
the node at $(k,M-L)$ representing $f(M,L,k)$. 
All incoming arrows at a given node indicate which 
values of $k$ and $M-L$ are needed in the recursion relations. For example
the nodes on the diagonal $k=M-L$ all follow from just the previous node
on the diagonal as given by (\ref{boundrec}). All other points need two
values of $k$ and $M-L$ as given by (\ref{rec}).

\vspace{5mm}
\proofarg{Proof of the initial condition (\ref{initial}).}{
\newline
For $k=0$ and $L=M$, $f_1$ and $f_2$ reduce to
\begin{equation}
f_1(M,M,0)=\sv{\vet}{M}
q^{\vet \, C^{-1}(\vet-\ve_{\la})}\prod_{j=1}^{N-1}
\Mults{\mu_j+\eta_j}{\eta_j},
\label{ia}
\end{equation}
with $\vmu=C^{-1}((2M+\ell)\ve_{N-1}+\ve_{\la}-2\vet)$, and
\begin{equation}
f_2(M,M,0)=\sv{\vn}{M}
q^{\vn \,  C^{-1}(\vn-\ve_{\la})}\prod_{j=1}^{N-1}\Mult{m_j+n_j}{n_j},
\end{equation}
with $\vm=C^{-1}((2M+\ell)\ve_{N-1}+\ve_{\la}-2\vn)$.
Since $\eta_j+\mu_j\geq 0$,
we can remove the prime in (\ref{ia}) and we are done.}
\vspace{5mm}

Since the $(\vm,\vn)$-system (\ref{mn}) corresponding to the
right-hand side of (\ref{lem}) does not depend on the variables $M-L$ and $k$,
$f_2$ is much easier to handle than $f_1$. Hence we proceed in reversed
order.

\proofarg{Proof of the recursion relation (\ref{rec}) for $f_2$}{
\newline
We substitute the definition of $f_2$, given by the
right-hand side of (\ref{lem}), into
$f_2(M,L,k)-(1-q^{M+L+\ell})f_2(M-1,L,k)$.
This immediately yields the desired result since
\begin{eqnarray}
\lefteqn{
f_2(M,L,k)-(1-q^{M+L+\ell})f_2(M-1,L,k) }\nonumber\\[2mm]
&&=
\sum_{r=L+k}^M q^{(r+L+\ell)(r-L-Nk)/N}
\left(\Mult{M-L-k}{M-r}-\Mult{M-L-k-1}{M-r-1}\right) \nonumber \\
& & \hspace{3cm} \times
\frac{(q^{\ell+1})_{L+M}}{(q^{\ell+1})_{L+r}}
\sv{\vn}{r}
q^{\vn \, C^{-1}(\vn-\ve_{\la})}\prod_{j=1}^{N-1}\Mult{m_j+n_j}{n_j}
\nonumber\\
&&=q^{M+L+\ell}\, f_2(M,L,k+1),
\end{eqnarray}
where in the last step we have used equation (\ref{arr}).}

\proofarg{Proof of the recursion relation (\ref{boundrec}) for $f_2$}{
\newline
The proof of (\ref{boundrec}) for $f_2$ is again trivial, since
\begin{eqnarray}
\lefteqn{q^{-(2L+\ell+1)(N-1)/N}\, f_2(M,L+1,M-L-1)} \nonumber\\
&& \qquad =q^{(M+L+\ell)(M-L-N(M-L))/N}
\sv{\vn}{M} q^{\vn \, C^{-1} (\vn-\ve_{\la})}
\prod_{j=1}^{N-1} \Mult{m_j+n_j}{n_j} \nonumber \\
&&\qquad = f_2(M,L,M-L).
\end{eqnarray}}

So far, everything has been extremely simple. The proof that
$f_1$ satisfies the recurrences (\ref{rec}) and (\ref{boundrec})
is however considerably more involved.

\proofarg{Proof of the recursion relations (\ref{rec}) for $f_1$}{
\newline
To avoid possible confusion in the course of the proof, we will
make the dependence on $M,L,k$ and $\vi$ of the vector $\vmu$ explicit
by writing $\vmu(M,L,k,\vi)$. Most importantly, the left-hand side of
(\ref{muetai}) should now be read as $\vmu(M,L,k,\vi)$.
We also make repeated use of the matrix elements of the
inverse Cartan matrix, given by
\begin{equation}
C^{-1}_{j,k} = \left\{\begin{array}{ll}
(N-k)j/N \qquad & {\rm for} \; j \leq k \\[2mm]
C^{-1}_{k,j} & {\rm for} \; j > k.
\end{array} \right.
\label{IC}
\end{equation}
\indent
{}From the definition of $f_1$ as the left-hand side of (\ref{lem})
we have
\begin{eqnarray}
\lefteqn{f_1(M,L,k)-f_1(M-1,L,k)} \nonumber\\[2mm]
&=& \svb{\vet,\vi}{L}
q^{-\vi \, C^{-1}(\vi+(2L+\ell)\ve_1)+\vet \,C^{-1}(\vet-\ve_{\la})}
\Mult{k}{\vi}_{\! 1/q}
\nonumber \\
&&\mbox{} \times
\Biggl\{ \prod_{j=1}^{N-1}\Mults{\eta_j+\mu_j(M,L,k,\vi)}{\eta_j}
-\prod_{j=1}^{N-1}\Mults{\eta_j+\mu_j(M,L,k,\vi)-1}{\eta_j}
\Biggr\}.
\label{rc}
\end{eqnarray}
Here we have used $\mu_j(M-1,L,k,\vi)=\mu_j(M,L,k,\vi)-C^{-1}_{1,j}-
C^{-1}_{N-1,j}=\mu_j(M,L,k,\vi)-1$.
\newline \indent
We now insert the telescopic expansion (\ref{rtele}) with
the replacements $A_j\to\eta_j$ and $B_j \to \mu_j$,
and change $\eta_p\to \eta_p+1$ and $i_p\to i_p-1$ in the
$p$-th term of the sum.
Carrying out both variable changes, using (\ref{muetai}) and
definition (\ref{iN}) of $i_N$, yields
\begin{eqnarray}
\lefteqn{
\mu_j(M,L,k,\vi) \to
\mu_j(M,L,k,\vi)
+ \ve_j C^{-1}(\ve_{p-1}-\ve_p-\ve_{N-1})} \nonumber \\
& & \hspace{3cm} = \mu_j(M,L,k,\vi)-\chi(j\geq p).
\end{eqnarray}
As a result (\ref{rc}) becomes (note that 
$\sum_{\vet+\ve_p,\vi-\ve_p}^{(\sigma)} \ldots =
\sum_{\vet,\vi}^{(\sigma)} \ldots$)
\begin{eqnarray}
\lefteqn{f_1(M,L,k)-f_1(M-1,L,k)} \nonumber\\
&=& \sum_{p=1}^{N-1} \svb{\vet,\vi}{L}
q^{-\vi\, C^{-1}(\vi+(2L+\ell)\ve_1)+\vet \, C^{-1} (\vet-\ve_{\la})}
\Mult{k}{\vi-\ve_p}_{\! 1/q} \nonumber \\[2mm]
&& \mbox{} \times
q^{\mu_p(M,L,k,\vi)-1+\ve_p C^{-1}
(2\vet +2\vi +(2L+\ell)\ve_1-\ve_{\la})}
\prod_{j=1}^{N-1} \Mults{\eta_j+\mu_j(M,L,k,\vi)-1}{\eta_j}.
\label{rd}
\end{eqnarray}
Using (\ref{muetai}) and (\ref{IC}) one may readily check that
\begin{eqnarray}
\mu_p(M,L,k,\vi) + \ve_p \, C^{-1} (2\vet + 2\vi +(2L+\ell)\ve_1 -\ve_{\la})
=M+L+\ell-k+\sum_{j=1}^p i_j,
\end{eqnarray}
resulting in
\begin{eqnarray}
\lefteqn{f_1(M,L,k)-f_1(M-1,L,k)}\nonumber\\
&=&q^{M+L+\ell}
\svb{\vet,\vi}{L}
q^{-\vi\, C^{-1}(\vi+(2L+\ell)\ve_1)}
\Biggl(\sum_{p=1}^{N-1} q^{-k-1+\sum_{j=1}^p i_j}
\Mult{k}{\vi-\ve_p}_{\! 1/q}\Biggr) \nonumber \\
&& \hspace{2cm} \times
q^{\vet \, C^{-1}(\vet-\ve_{\la})}
\prod_{j=1}^{N-1} \Mults{\eta_j+\mu_j(M,L,k,\vi)-1}{\eta_j}.
\label{rg}
\end{eqnarray}
\indent
One can now simplify the sum over $p$. Using definition (\ref{iN})
we get
\begin{eqnarray}
\sum_{p=1}^{N-1} q^{k+1-\sum_{j=1}^p i_j}
\Mult{k}{\vi-\ve_p}_q 
& = & \frac{q^{k+1}(q)_k}{(q)_{i_1}\ldots
(q)_{i_{N-1}}(q)_{i_N+1}}
\sum_{p=1}^{N-1} (1-q^{i_p}) q^{-\sum_{j=1}^p i_j} \nonumber\\
&=& \frac{(q)_k}{(q)_{i_1}\ldots
(q)_{i_{N-1}}(q)_{i_N+1}}
\Bigl((1-q^{k+1})-(1-q^{i_N+1})\Bigr) \nonumber \\[2mm] 
&=&\Mult{k+1}{\vi}_q -\Mult{k}{\vi}_q .
\end{eqnarray}
Substituting this into (\ref{rg}) and noting
$\mu(M,L,k,\vi)-1=\mu(M,L,k+1,\vi)=\mu(M-1,L,k,\vi)$, we arive at
the desired result
\begin{eqnarray}
f_1(M,L,k)-f_1(M-1,L,k) =
q^{M+L+\ell}\Bigl(f_1(M,L,k+1)-f_1(M-1,L,k) \Bigr).
\end{eqnarray}}

\proofarg{Proof of recursion relation (\ref{boundrec}) for $f_1$}{\newline 
{}From the definition of $f_1$ as the left-hand side of (\ref{lem})
we obtain
\begin{eqnarray}
\lefteqn{f_1(M,L+1,M-L-1) =
\svb{\vet,\vi}{L+1}
q^{-\vi \, C^{-1}(\vi+(2L+\ell+2)\ve_1)}
\Mult{M-L-1}{\vi}_{\! 1/q}}\nonumber\\
\lefteqn{ \hspace{3cm}
\times
q^{ \vet \, C^{-1} (\vet-\ve_{\la})}
\prod_{j=1}^{N-1} \Mults{\eta_j+\mu_j(M,L+1,M-L-1,\vi)}{\eta_j}.}
\end{eqnarray}
We now apply telescopic expansion (\ref{btele}),
with $A_j \to \eta_j$ and $B_j \to \mu_j$,
to rewrite the product
over the $q$-binomials. 
Then we make the change of  
variables $i_1\to i_1-1$ in the first-term on the
right-hand side of (\ref{btele}),
$\eta_p\to \eta_p+1,\;i_{p+1}\to
i_{p+1}-1$ in the $p$-th term of the sum $(1\leq p\leq N-2)$
and $\eta_{N-1} \to \eta_{N-1}+1$ in the term with $p=N-1$.
Since under these variable changes
\begin{eqnarray}
\lefteqn{
\mu_j(M,L+1,M-L-1,\vi) \to
\mu_j(M,L,M-L,\vi)
- \ve_j C^{-1}(\ve_p-\ve_{p+1}-\ve_{N-1})} \nonumber \\
& & \hspace{5cm} = \mu_j(M,L,M-L,\vi)+\chi(j>p),
\end{eqnarray}
we obtain
\begin{eqnarray}
\lefteqn{f_1(M,L+1,M-L-1)} \nonumber \\
&=&
\sum_{p=0}^{N-1}
\svb{\vet,\vi}{L}
q^{-\vi \, C^{-1}(\vi+(2L+\ell)\ve_1) +
\vet \, C^{-1} (\vet-\ve_{\la})}
\Mult{M-L-1}{\vi-\ve_{p+1}}_{\! 1/q} \nonumber\\[2mm]
& & \!\! \times \, q^{\mu_p(M,L,M-L,\vi)
-2\vi C^{-1} (\ve_1-\ve_{p+1})
+\ve_p C^{-1} (2\vet-\ve_{\la}+\ve_p)
+\ve_{p+1}C^{-1} ((2L+\ell+2)\ve_1-\ve_{p+1}) }
\nonumber\\[2mm]
&& \! \! \times
\prod_{j=1}^{N-1} \Mults{\eta_j+\mu_j(M,L,M-L,\vi)}{\eta_j}.
\label{bd}
\end{eqnarray}
Now use (\ref{muetai}) as well as the definition of $i_N$ (recall that $k=M-L$)
to observe that the exponent of $q$ in the
third line of the above equation simplifies to
\begin{equation}
(2L+\ell+1)(N-1)/N -\sum_{j=1}^p i_j.
\end{equation}
This allows (\ref{bd}) to be reduced to
\begin{eqnarray}
\lefteqn{q^{-(2L+\ell+1)(N-1)/N} f_1(M,L+1,M-L-1)} \nonumber\\
&=& 
\svb{\vet,\vi}{L}
q^{-\vi\, C^{-1}(\vi+(2L+\ell)\ve_1)}
\: \Biggl( \: \sum_{p=0}^{N-1} q^{-\sum_{j=1}^p i_j}
\Mult{M-L-1}{\vi-\ve_{p+1}}_{1/q}\Biggr) \nonumber\\
&& \times 
q^{\vet \, C^{-1}(\vet-\ve_{\la})}
\Mults{\eta_j+\mu_j(M,L,M-L,\vi)}{\eta_j}.
\label{bf}
\end{eqnarray}
Again we may simplify the sum over $p$,
\begin{eqnarray}
\sum_{p=0}^{N-1} q^{\sum_{j=1}^p i_j}
\Mult{M-L-1}{\vi-\ve_{p+1}}_q &=& \frac{(q)_{M-L-1}}{(q)_{i_1}\ldots
(q)_{i_N}} \sum_{p=0}^{N-1}
(1-q^{i_{p+1}}) q^{\sum_{j=1}^p i_j} \nonumber\\[2mm]
&=&\Mult{M-L}{\vi}_q,
\end{eqnarray}
so that the right-hand side of (\ref{bf}) gives $f_1(M,L,M-L)$ as desired.}

\nsection{Applications}\label{secapp}
We now show how theorem~\ref{thGD}, and in particular
corollary~\ref{corgd}, can be exploited to obtain very
general $q$-series identities, which take many of the
known Rogers--Ramanujan type identities as special cases.

We substitute the conjugate Bailey pair
of corollary~\ref{corgd} into equation~\ref{agbd}
and take the limit $M\to\infty$.
This gives the important result
\begin{corollary}[higher-level Bailey lemma]\label{corinf}
Let $N,\ell,\la$ and $\sigma$ be as in theorem~\ref{thGD},
and let $(\alpha,\beta)$ form a Bailey pair relative to $a=q^{\ell}$.
Then the following identity holds:
\begin{eqnarray}
\lefteqn{
\frac{1}{(a q)_{\infty}}
\sum_{L=0}^{\infty}
a^{L/N} q^{L^2/N} \, \alpha_L
\sv{\vet}{L}
\frac{q^{\vet \, C^{-1} (\vet-\ve_{\la})}}
{(q)_{\eta_1} \ldots (q)_{\eta_{N-1}}} } \nonumber \\
&=&
\sum_{L=0}^{\infty}
a^{L/N} q^{L^2/N} \, \beta_L
\sv{\vn}{L}
q^{\vn \, C^{-1} (\vn-\ve_{\la})}
\prod_{j=1}^{N-1} \Mult{m_j+n_j}{n_j},
\label{eqinf}
\end{eqnarray}
with the $(\vm,\vn)$-system
\begin{equation}
\vm = C^{-1} \Bigl((2L+\ell)\,\ve_{N-1}+\ve_{\la}-2\vn\Bigr).
\end{equation}
\end{corollary}
For $\la=(\ell),~\sigma=0$ and for $|\la|=0,~\ell=N,~\sigma=1$ this
result was announced in ref.~\cite{SW}.

As first noted by Foda and Quano~\cite{FQ} for some specific examples and in
full in refs.~\cite{BM,BMS96}, polynomial
identites related to the characters of $c<1$ Virasoro 
algebras give rise to Bailey pairs.
The name of the above lemma has been motivated by the fact that
inserting these Bailey pairs into (\ref{eqinf}) 
yields Rogers--Ramanujan type identities for
the branching functions of level-$N$ A$_1^{(1)}$ conformal
coset models.
Details of this application 
will be presented elsewhere~\cite{BMSW}.\footnote{It is intriguing
to note~\cite{BMS95} that for $N=2$ the characters of these coset models
(which are obtained from the branching functions by summing over 
both parities of $\sigma$) also arise from the ``classical'' conjugate 
Bailey pair (\ref{gammarho}) and (\ref{deltarho}), taking $\rho_1, M
\to \infty$ and setting $\rho_2=-(aq)^{1/2}$.} 
Instead, we apply here corollary \ref{corinf} to some of the simplest
known Bailey pairs. For these we wish to use corollary \ref{corinf}
as efficiently as possible.
To this end we utilize, as well as the Bailey chain
and its generalization to the Bailey lattice, two new transformations
acting on Bailey pairs.

\subsection{Transformations on Bailey pairs}
Lemma~\ref{LAB} provides a powerful mechanism to
obtain new Bailey pairs from a given Bailey pair.
One of the limitations of (\ref{AB}) is however that if
\bp \ is a Bailey pair relative to $a$, so will be
the resulting Bailey pair \bpp. 

A useful transformation which does not keep $a$ fixed 
was found by Agarwal, Andrews and Bressoud~\cite{AAB,Bressoud88}
\begin{lemma}\label{lemBlattice}
Let $(\alpha,\beta)$ be a Bailey pair relative to $a$
and define $\alpha'_0=\alpha_0$ and 
\begin{eqnarray}
\alpha'_L &=&  
\frac{(\rho_1)_L (\rho_2)_L (a/\rho_1 \rho_2)^L}
{(a/\rho_1)_L (a/\rho_2)_L}
\left( \frac{(a)_{2L}}{(aq)_{2L}} \: \alpha_L
-a q^{2L-2} \frac{(a)_{2L-2}}{(aq)_{2L-2}} \: \alpha_{L-1} \right)
\qquad  L \geq 1
\nonumber \\[1mm]
\beta'_L &=& \sum_{r=0}^L
\frac{(\rho_1)_r (\rho_2)_r (a/\rho_1 \rho_2)^r
(a/\rho_1 \rho_2)_{L-r}}
{(a/\rho_1)_L (a/\rho_2)_L (q)_{L-r}} \: \beta_r.
\label{EqBl}
\end{eqnarray}
Then $(\alpha',\beta')$ is a Bailey pair relative to $aq^{-1}$.
\end{lemma}
In combination with lemma~\ref{LAB}, transformation (\ref{EqBl}) gives
rise to what is called the Bailey lattice.

We now introduce two new transformations that map a Bailey pair onto
a new Bailey pair.
The first transformation is of the type (\ref{AB}) in that it leaves
the value of $a$ invariant.
\begin{lemma}\label{lemchain}
If $(\alpha,\beta)$ forms a Bailey pair relative to $a$, then
the sequences
\begin{eqnarray}
\alpha'_L &=&   a^L q^{L(L+1)}
\sum_{r=0}^L
\frac{(\rho_1)_r (\rho_2)_r (q^{1-r}/\rho_1 \rho_2)^r}
{(aq/\rho_1)_r (aq/\rho_2)_r}
\: \alpha_r  \nonumber \\
&&-a^{L-1} q^{L(L-1)}
\sum_{r=0}^{L-1}
\frac{(\rho_1)_r (\rho_2)_r (q^{1-r}/\rho_1 \rho_2)^r}
{(aq/\rho_1)_r (aq/\rho_2)_r}
\: \alpha_r  \nonumber \\
\beta'_L &=& q^L \sum_{r=0}^L
\frac{(\rho_1)_r (\rho_2)_r (aq/\rho_1 \rho_2)^r
(aq/\rho_1 \rho_2)_{L-r}}
{(aq/\rho_1)_L (aq/\rho_2)_L (q)_{L-r}} \: \beta_r .
\label{Eqchain2}
\end{eqnarray}
again form a Bailey pair relative to $a$.
\end{lemma}
The second transformation closely mimics lemma~\ref{lemBlattice}.
\begin{lemma}\label{lemBl2}
Let $(\alpha,\beta)$ be a Bailey pair relative to $a$
and define $\alpha'_0=\alpha_0$ and
\begin{eqnarray}
\alpha'_L &=&  
\frac{(\rho_1)_L (\rho_2)_L (aq/\rho_1 \rho_2)^L}
{(aq/\rho_1)_L (aq/\rho_2)_L}
\frac{(a)_{2L}}{(aq)_{2L}} \: \alpha_L  \nonumber \\[1mm]
&& - \frac{(\rho_1)_{L-1} (\rho_2)_{L-1} (aq/\rho_1 \rho_2)^{L-1}}
{(aq/\rho_1)_{L-1} (aq/\rho_2)_{L-1}} \: aq^{2L-2}
\frac{(a)_{2L-2}}{(aq)_{2L-2}} \: \alpha_{L-1}
\qquad  L \geq 1 \nonumber \\
\beta'_L &=& \sum_{r=0}^L
\frac{(\rho_1)_r (\rho_2)_r (aq/\rho_1 \rho_2)^r
(aq/\rho_1 \rho_2)_{L-r}}
{(aq/\rho_1)_L (aq/\rho_2)_L (q)_{L-r}} \: \beta_r .
\label{EqBl2}
\end{eqnarray}
Then $(\alpha',\beta')$ is a Bailey pair relative to $aq^{-1}$.
\end{lemma}

\proofarg{Proof of the lemmas~\ref{lemchain} and \ref{lemBl2}}{\newline
For both the proof of equation (\ref{Eqchain2}) and (\ref{EqBl2})
we only need the
$q$-analogue of Saalsch\"utz's theorem 
(equation (3.3.12) of ref.~\cite{AndrewsBook}).
First, to prove (\ref{Eqchain2}),
\begin{eqnarray}
\beta'_L &=& q^L \sum_{k=0}^L
\frac{(\rho_1)_k (\rho_2)_k (aq/\rho_1 \rho_2)^k
(aq/\rho_1 \rho_2)_{L-k}}
{(aq/\rho_1)_L (aq/\rho_2)_L (q)_{L-k}} \: \beta_k  \nonumber \\
&=& q^L \sum_{k=0}^L
\frac{(\rho_1)_k (\rho_2)_k (aq/\rho_1 \rho_2)^k
(aq/\rho_1 \rho_2)_{L-k}}
{(aq/\rho_1)_L (aq/\rho_2)_L (q)_{L-k}}
\sum_{m=0}^k \frac{\alpha_m}{(q)_{k-m}(aq)_{k+m}} \nonumber \\
&=& q^L\sum_{m=0}^L \alpha_m
\sum_{k=m}^L
\frac{(\rho_1)_k (\rho_2)_k (aq/\rho_1 \rho_2)^k
(aq/\rho_1 \rho_2)_{L-k}}
{(aq/\rho_1)_L (aq/\rho_2)_L (q)_{L-k}(q)_{k-m}(aq)_{k+m}}
\nonumber \\
&=& q^L \sum_{m=0}^L
\frac{(\rho_1)_m (\rho_2)_m (aq/\rho_1 \rho_2)^m \, \alpha_m}
{(aq/\rho_1)_m (aq/\rho_2)_m (q)_{L-m}(aq)_{L+m}} 
\nonumber \\
&=& \sum_{m=0}^L
\frac{(\rho_1)_m (\rho_2)_m (q^{1-m}/\rho_1 \rho_2)^m \, \alpha_m}
{(aq/\rho_1)_m (aq/\rho_2)_m} 
\Biggl( \sum_{k=m}^L aq^{2k} - \sum_{k=m+1}^L \Biggr)
\frac{a^{k-1} q^{k(k-1)}}{(q)_{L-k}(aq)_{L+k}} 
\nonumber \\
&=&
\sum_{k=0}^L
\frac{a^k q^{k(k+1)}}{(q)_{L-k}(aq)_{L+k}}
\sum_{m=0}^k
\frac{(\rho_1)_m (\rho_2)_m (q^{1-m}/\rho_1 \rho_2)^m \, \alpha_m}
{(aq/\rho_1)_m (aq/\rho_2)_m}
\nonumber \\
&& -
\sum_{k=0}^L
\frac{a^{k-1} q^{k(k-1)}}{(q)_{L-k}(aq)_{L+k}}
\sum_{m=0}^{k-1}
\frac{(\rho_1)_m (\rho_2)_m (q^{1-m}/\rho_1 \rho_2)^m \, \alpha_m}
{(aq/\rho_1)_m (aq/\rho_2)_m}
\nonumber \\
&=& \sum_{k=0}^L \frac{\alpha_k'}
{(q)_{L-k}(aq)_{L+k}} \, .
\end{eqnarray}
To prove (\ref{EqBl2}) we proceed as follows,
\begin{eqnarray}
\beta'_L &=&
\sum_{k=0}^L
\frac{(\rho_1)_k (\rho_2)_k (aq/\rho_1 \rho_2)^k
(aq/\rho_1 \rho_2)_{L-k}}
{(aq/\rho_1)_L (aq/\rho_2)_L (q)_{L-k}} \: \beta_k  \nonumber \\
&=& \sum_{m=0}^L
\frac{(\rho_1)_m (\rho_2)_m (aq/\rho_1 \rho_2)^m \, \alpha_m}
{(aq/\rho_1)_m (aq/\rho_2)_m (q)_{L-m}(aq)_{L+m}} 
\nonumber \\
&=& \sum_{m=0}^L
\frac{(\rho_1)_m (\rho_2)_m (aq/\rho_1 \rho_2)^m (a)_{2m} \, \alpha_m}
{(aq/\rho_1)_m (aq/\rho_2)_m (aq)_{2m}}
\Biggl( \frac{1}{(q)_{L-m}(a)_{L+m}}-
\frac{aq^{2m}}{(q)_{L-m-1}(a)_{L+m+1}}\Biggr)
\nonumber \\
&=& \sum_{m=0}^L \frac{\alpha_m'}
{(q)_{L-m}(a)_{L+m}} \, .
\end{eqnarray}}

\subsection{A general $q$-series identity}
Equipped with corollary~\ref{corinf} and the lemmas
\ref{LAB}, \ref{lemBlattice}--\ref{lemBl2},
we are prepared to prove the following $q$-identity.
\begin{theorem}\label{thidentity}
Fix integers $N \geq 1$, $\delta \in \{0,1\}$, $k\geq 2$ and $1\leq i \leq k$.
Let $\la$ be a partition with largest part $\leq N-1$
and $\sigma$ $\in \{0,1\}$ such that $|\la|+\sigma N$ even.
Then the following identity holds:
\begin{eqnarray}
\lefteqn{
\frac{1}{(q)_{\infty}}
\sum_{j=-\infty}^{\infty}
(-1)^j
q^{\bigl((2k+\delta-2+2/N)j+2k-2i+\delta\bigr)j/2}
\sv{\vet}{j} \frac{q^{\vet\, C^{-1} (\vet-\ve_{\la})}}
{(q)_{\eta_1} \ldots (q)_{\eta_{N-1}}} } \nonumber \\[2mm]
&=&
\sum_{r_1 \geq \ldots  \geq r_{k-1} \geq 0}
\frac{
q^{r_1^2/N + r_2^2 + \cdots + r_{k-1}^2 + r_i+\cdots + r_{k-1}}}
{(q)_{r_1-r_2} \ldots (q)_{r_{k-2}-r_{k-1}}
(q^{2-\delta};q^{2-\delta})_{r_{k-1}} }\nonumber \\[1mm]
& & \hspace{3cm} \times
\sv{\vn}{r_1} q^{\vn\, C^{-1}(\vn-\ve_{\la})}
\prod_{j=1}^{N-1} \Mult{m_j+n_j}{n_j} 
\label{Eqq}
\end{eqnarray}
with $(\vm,\vn)$-system $ \vm = C^{-1} (2r_1\,\ve_{N-1}+\ve_{\la}-2\vn)$.
\end{theorem}
(Note: the summation symbol $\sum^{(\sigma,a)}$ in the above is that
defined in (\ref{sumsymbol}) with $\ell=0$.)

Before we prove this result, we note that identity (\ref{Eqq})
has a number of important specializations.
First, setting $N=1$ and rewriting the left-hand side
using Jacobi's triple product identity 
(equation (2.2.10) of ref.~\cite{AndrewsBook}), we get
\begin{corollary}
For integers 
$\delta=0,1$, $k\geq 2$ and $1\leq i \leq k+\delta-1$,
\begin{equation}
\sum_{n_1,\ldots,n_{k-1} \geq 0}
\frac{
q^{N_1^2 + \cdots + N_{k-1}^2 + N_i+\cdots + N_{k-1}}}
{(q)_{n_1} \ldots (q)_{n_{k-2}}
(q^{2-\delta};q^{2-\delta})_{n_{k-1}} }
\hspace{4mm}= \hspace{-4mm} \prod_{\begin{array}{c}
\\[-6mm]
\scriptstyle j=1 \\[-1mm]
\scriptstyle j \not\equiv 0, \pm i \mod{2k+\delta}
\end{array}}^{\infty} \hspace{-12mm} (1-q^j)^{-1},
\label{N1}
\end{equation}
where $N_j = n_j + \cdots + n_{k-1}$.
\end{corollary}
For $\delta=1$, $k=2$ these are the Rogers--Ramanujan 
identities~\cite{Rogers,RR}.
For $\delta=1$, general $k$, these are Andrews' analytic counterpart
of Gordon's partition theorem~\cite{Andrews74,Gordon61}.
For $\delta=0$ the identities (\ref{N1}) were discovered by
Bressoud~\cite{Bressoud80},
and are the analytic counterpart of 
partition identities of Andrews~\cite{Andrews67}
and Bressoud~\cite{Bressoud79}.
We note that the derivation of these identities using
the Bailey chain and lattice is of course not new,
see e.g., ref.~\cite{AAB}.

The next two corollaries are obtained by 
setting $N=2$ in theorem~\ref{thidentity}
\begin{corollary}\label{corN2a}
For integers $k\geq 2$ and $1\leq i \leq k$,
\begin{eqnarray}
\sum_{n_1,\ldots,n_{k-1}\geq 0}
\frac{q^{N_1^2 + 2N_2^2 + \cdots + 2N_{k-1}^2 + 2N_i+\cdots + 2N_{k-1}}
(-q;q^2)_{N_1}}
{(q^2;q^2)_{n_1} \ldots (q^2;q^2)_{n_{k-1}}}
\hspace{4mm}= \hspace{-8mm} \prod_{\begin{array}{c}
\\[-6mm]
\scriptstyle j=1 \\[-1mm]
\scriptstyle j \not\equiv 2 \mod{4} \\[-1mm]
\scriptstyle j \not\equiv 0, \pm (2i-1) \mod{4k}
\end{array}}^{\infty} \hspace{-12mm}
(1-q^j)^{-1}.
\label{N2a}
\end{eqnarray}
\end{corollary}
\begin{corollary}\label{corN2b}
For integers 
$k\geq 2$ and $1\leq i \leq k-1$,
\begin{eqnarray}
\lefteqn{
\sum_{n_1,\ldots,n_{k-1}\geq 0}
\frac{q^{N_1^2 + 2N_2^2 + \cdots + 2N_{k-1}^2 + 2N_i+\cdots + 2N_{k-1}}
(-q;q^2)_{N_1}}{(q^2;q^2)_{n_1} \ldots (q^2;q^2)_{n_{k-2}}
(q^4;q^4)_{n_{k-1}}} } \nonumber \\
&& \hspace{4cm} = 
(-q^{2k-1},q^{4k-2})_{\infty} \hspace{-5mm}
\hspace{-5mm} \prod_{\begin{array}{c}
\\[-6mm]
\scriptstyle j=1 \\[-1mm]
\scriptstyle j \not\equiv 2 \mod{4} \\[-1mm]
\scriptstyle j \not\equiv 0 \mod{8k-4} \\[-1mm]
\scriptstyle j \not\equiv 2k-1, \pm (2i-1) \mod{4k-2}
\end{array}}^{\infty} \hspace{-16mm}
(1-q^j)^{-1}.
\label{N2b}
\end{eqnarray}
\end{corollary}
The identities of corollary~\ref{corN2a} are due to Andrews~\cite{Andrews75} 
and Bressoud~\cite{Bressoud80b},
and are related to Andrews' generalization of the G\"ollnitz--Gordon
partition identities~\cite{Goellnitz,Gordon65,Andrews67b}.
The identities of corollary~\ref{corN2b} were discovered 
by Bressoud~\cite{Bressoud80b}.
These results can also be obtained using the conjugate Bailey pair given 
by~(\ref{gammarho}) and~(\ref{deltarho}). 
For $N\geq 3$, however, theorem~\ref{thidentity} relies on the higher-level 
Bailey lemma and, at present, we cannot obtain it by means of the standard 
conjugate Bailey pair.

\proofarg{Proof of corollaries~\ref{corN2a} and \ref{corN2b}}{\newline
Set $N=2$ in (\ref{Eqq}) and sum over both choices for $\sigma$.
Recalling the $q$-binomial expansion
(equation (3.3.6) of ref.~\cite{AndrewsBook}) to sum
\begin{equation}
\sum_{n=0}^{r_1+|\la|/2}
q^{n(n-|\la|)/2} \Mult{r_1+|\la|/2}{n} = 
(-q^{(1-|\la|)/2})_{|\la|/2}(-q^{1/2})_{r_1}
\end{equation}
and
\begin{equation}
\sum_{\eta=0}^{\infty}
\frac{q^{\eta(\eta-|\la|)/2}}{(q)_{\eta}} =
(-q^{(1-|\la|)/2})_{|\la|/2}(-q^{1/2})_{\infty}
\end{equation}
we get, after dropping the factor 
$(-q^{(1-|\la|)/2})_{|\la|/2}$,
replacing $q\to q^2$ and applying the
triple product identity,
\begin{eqnarray}
\lefteqn{
\sum_{n_1,\ldots,n_{k-1}\geq 0}
\frac{q^{N_1^2 + 2N_2^2 + \cdots + 2N_{k-1}^2 + 2N_i+\cdots + 2N_{k-1}}
(-q;q^2)_{N_1}}{(q^2;q^2)_{n_1} \ldots (q^2;q^2)_{n_{k-2}}
(q^{4-2\delta};q^{4-2\delta})_{n_{k-1}}} } \nonumber \\
&=& \frac{1}{(q)_{\infty}}
\prod_{n=0}^{\infty}
(1-q^{2+4n})
(1-q^{2i-1+An})
(1-q^{-2i+1+A(n+1)})
(1-q^{A(n+1)}),
\end{eqnarray}
with $A=4k+2\delta-2$. Rewriting the right-hand side gives
equation (\ref{N2a}) when $\delta=1$ and equation (\ref{N2b}) when
$\delta=0$.}

\subsection{Proof of theorem~\ref{thidentity}}
The proof of theorem~\ref{thidentity} is based on the following
two Bailey pairs relative to $q$
(equations (2.12) and (2.13) of ref.~\cite{Andrews84} and
item E(3) of ref.~\cite{Slater}):
\begin{eqnarray}
\alpha_L &=& (-1)^L q^{L(L-1)/2} \:
\frac{(q^2)_{2L}}{(q)_{2L}} \nonumber \\[2mm]
\beta_L &=& \delta_{L,0}
\label{I}
\end{eqnarray}
and
\begin{eqnarray}
\alpha_L &=& (-1)^L q^{L^2} \: 
\frac{(q^2)_{2L}}{(q)_{2L}} \nonumber \\[2mm]
\beta_L &=& \frac{1}{(q^2;q^2)_L},
\label{II}
\end{eqnarray}
as well as the Bailey pair (items B(2) and E(4) of ref.~\cite{Slater}) 
\begin{eqnarray}
\alpha_L &=& \left\{\begin{array}{ll}
1 &  L=0 \\[2mm]
(-1)^L q^{(\delta+2)L(L-1)/2} \Bigl(1+q^{(\delta+2)L}\Bigr) 
\quad & L\geq 1
\end{array} \right. \nonumber \\
\beta_L &=& \frac{q^L}{(q^{2-\delta};q^{2-\delta})_L},
\label{III}
\end{eqnarray}
relative to 1, where $\delta=0,1$. 

To now derive (\ref{Eqq}) we have to consider various 
different ranges of $i$.

\proofarg{Proof of theorem~\ref{thidentity} for $i=1$.}{\newline
We take lemma~\ref{lemchain} with $a=1$ and let $\rho_1,\rho_2\to\infty$.
Hence given a Bailey pair \bp \ relative to 1, the sequences
\begin{eqnarray}
\alpha'_L &=&   
q^{L(L+1)}\alpha_L
-q^{L(L-1)}(1-q^{2L})(\alpha_{L-1} + \cdots + \alpha_0)
\nonumber \\[2mm]
\beta'_L &=& q^L \sum_{k=0}^L
\frac{q^{k^2} \, \beta_k}{(q)_{L-k}}
\end{eqnarray}
form yet again a Bailey pair relative to 1.
Now iterating this $k-2$ times, using (\ref{III})
as initial Bailey pair, we obtain the new Bailey pair
\begin{eqnarray}
\alpha_L &=& \left\{\begin{array}{ll}
1 &  L=0 \\[2mm]
(-1)^L q^{(2k+\delta-2)L(L-1)/2} \Bigl(1+q^{(2k+\delta-2)L}\Bigr)
\quad & L\geq 1
\end{array} \right. \nonumber \\[2mm]
\beta_L &=& q^L \sum_{L\geq r_2 \geq \ldots \geq r_{k-1}}
\frac{q^{r_2^2 + \cdots + r_{k-1}^2+r_2 + \cdots + r_{k-1}}}
{(q)_{L-r_2} (q)_{r_2-r_3} \ldots (q)_{r_{k-2}-r_{k-1}}
(q^{2-\delta};q^{2-\delta})_{r_{k-1}}} \, .
\label{nbp}
\end{eqnarray}
Substitution into equation~(\ref{eqinf}) with $a=1$
yields theorem~\ref{thidentity} for $i=1$.}

\proofarg{Proof of theorem~\ref{thidentity} for $2\leq i \leq k+\delta-1$.}
{\newline
We take the Bailey pair (\ref{I}) when $\delta=1$
and (\ref{II}) when $\delta=0$.
Now we apply transformation (\ref{AB}) $k-i+\delta-1$ times, 
then (\ref{EqBl2}) once and then (\ref{AB}) $i-2$ times, all
with $\rho_1,\rho_2 \to \infty$.
This yields the following Bailey pair relative to 1:
\begin{eqnarray}
\alpha_L &=& \left\{\begin{array}{ll}
1 &  L=0 \\[2mm]
(-1)^L q^{\bigl((2k+\delta-2)L-2k+2i-\delta\bigr)L/2} 
\Bigl(1+q^{(2k-2i+\delta)L}\Bigr)
\quad & L\geq 1
\end{array} \right. \nonumber \\[2mm]
\beta_L &=& \sum_{L\geq r_2 \geq \ldots \geq r_{k-1}}
\frac{q^{r_2^2 + \cdots + r_{k-1}^2+r_i + \cdots + r_{k-1}}}
{(q)_{L-r_2} (q)_{r_2-r_3} \ldots (q)_{r_{k-2}-r_{k-1}}
(q^{2-\delta};q^{2-\delta})_{r_{k-1}}} \, ,
\label{BPi}
\end{eqnarray}
for $i=2,\ldots,k+\delta-1$.
Substituting this into equation~(\ref{eqinf}) with $a=1$
completes the proof.}

\proofarg{Proof of theorem~\ref{thidentity} for $3\leq i \leq k$.}
{\newline
Again we take the Bailey pair (\ref{I}) when $\delta=1$
and (\ref{II}) when $\delta=0$.
This time we apply transformation (\ref{AB}) $k-i+\delta$ times,
then (\ref{EqBl}) once and then (\ref{AB}) $i-3$ times, all
with $\rho_1,\rho_2 \to \infty$.
This again yields the Bailey pair (\ref{BPi})
but now for $i=3,\ldots,k+\delta$.
Since the cases $i=k$ and $i=k+\delta$ coincide, we only
need to keep $i=3,\ldots,k$.
After substituting this into equation~(\ref{eqinf}) with $a=1$
we are done.}

\nsection{Conclusion}
We end this paper with a few further remarks about the
higher-level Bailey lemma.

In this paper we considered the problem of generalizing
the conjugate Bailey pair (\ref{gammarho})--(\ref{deltarho}),
and more specifically its specialization (\ref{gammadelta}).
Now the very problem of what exactly is meant by ``generalization''
deserves some discussion.
In particular we note the following lemma.
\begin{lemma}
Let \bp \ and \bpp \ be Bailey pairs relative to $a$ and $ab$,
respectively, with
\begin{equation}
\alpha'_L = \sum_{k=0}^L P_{L,k} \alpha_k \quad {\rm and} \quad
\beta'_L = \sum_{k=0}^L Q_{L,k} \beta_k,
\label{PQ}
\end{equation}
and let \gd \ be a conjugate Bailey pair relative to $a$.
Then the sequences
\begin{equation}
\gamma'_L = \sum_{k=L}^{\infty} P_{k,L} \gamma_k \quad {\rm and} \quad
\delta'_L = \sum_{k=L}^{\infty} Q_{k,L} \delta_k,
\label{PQ2}
\end{equation}
form a conjugate Bailey pair relative to $ab^{-1}$.
\end{lemma}
\proof{For compactness, write the equations~(\ref{ab}) and
(\ref{gd})
as $\beta_L = \sum_{k=0}^L M_{L,k}(a) \alpha_k$ and
$\gamma_L = \sum_{k=L}^{\infty} M_{k,L}(a) \delta_k$.
Then from the fact that \bp \ and \bpp \ are Bailey pairs, and from
equation (\ref{PQ}), we infer the relation
\begin{equation}
\sum_{k=m}^L Q_{L,k} M_{k,m}(a) = 
\sum_{k=m}^L M_{L,k}(ab) P_{k,m} \qquad m=0,\ldots,L.
\end{equation}
With this we compute
\begin{eqnarray}
\lefteqn{
\gamma_L' = \sum_{k=L}^{\infty} P_{k,L} \gamma_k = 
\sum_{k=L}^{\infty} P_{k,L} \sum_{m=k}^{\infty} M_{m,k}(a) \delta_m }
\nonumber \\
\lefteqn{ \hphantom{\gamma_L'} = 
\sum_{m=L}^{\infty} \delta_m \sum_{k=L}^m
M_{m,k}(a) P_{k,L}  =
\sum_{m=L}^{\infty} \delta_m \sum_{k=L}^m
Q_{m,k} M_{k,L}(a/b)  } \nonumber \\
\lefteqn{ \hphantom{\gamma_L'} = 
\sum_{k=L}^{\infty} M_{k,L}(a/b)
\sum_{m=k}^{\infty} Q_{m,k} \delta_m =
\sum_{k=L}^{\infty} M_{k,L}(a/b) \delta_k',}
\end{eqnarray}
establishing the claim of the lemma.}

The obvious implication of this result is that instead of
using the Bailey-chain, Bailey-lattice, et cetera, we could 
equally well keep \bp \ fixed, and transform \gd \ to yield
a conjugate Bailey-chain, lattice and so on.
Of course, when searching for generalizations of the
conjugate Bailey pair (\ref{gammadelta}), one must exclude
those conjugate Bailey pairs which can simply be obtained from 
(\ref{gammadelta}) by linear transformations of type
(\ref{PQ2}).
Let us merely state here that the new conjugate Bailey pairs
presented in this paper admit no further simplifications,
and are, in the loose sense of
ref.~\cite{Andrews84}, ``reduced''.

\vspace{5mm}
The expression on the left-hand side of the higher-level Bailey lemma
admits several different representations from the one listed in 
equation (\ref{eqinf}).
The most relevant rewriting arises by first noting that the 
sum over $\vet$ only depends on $L \mod{N}$.
Using this as well as some other simple symmetries, 
the left-hand side of equation~(\ref{eqinf})
can be put into the form
\begin{equation}
\frac{1}{(a q)_{\infty}} \sum_{p-p_{\ell}=0}^{\lfloor N/2 \rfloor}
\sum_{\vet \in \Integer^{N-1} \atop
\frac{2p-|\la|+\sigma N}{2N}-(C^{-1} \vet)_1 \in \Integer}
\hspace{-4mm} \frac{q^{\vet \, C^{-1} (\vet-\ve_{\la})}}
{(q)_{\eta_1} \ldots (q)_{\eta_{N-1}}}
\hspace{-5mm}
\sum_{\begin{array}{c}
\\[-6mm]
\scriptstyle L=0 \\[-1mm]
\scriptstyle L+\ell/2 \equiv \pm p \mod{N}
\end{array}}^{\infty} \hspace{-7mm}
a^{L/N} q^{L^2/N} \, \alpha_L,
\label{e55}
\end{equation}
where $p_{\ell} = 0$ when $\ell$ is even and
$p_{\ell} = 1/2$ when $\ell$ is odd.
We now assume that the partition $\la$ has just a single part,
$\la=(\ell')$, with $\ell'=0,\dots,N-1$.
This admits (\ref{e55}) to be expressed in terms of a sum over
the level-$N$ A$^{(1)}_1$ string functions $c_m^{\ell}$~\cite{JM,KP,LP}, 
\begin{equation}
q^{-\frac{\ell'(N-\ell')}{2N(N+2)}} 
(q)_{\ell} 
\sum_{p-p_{\ell}=0}^{\lfloor N/2 \rfloor}
c_{2p+\sigma N}^{\ell'}
\hspace{-5mm}
\sum_{\begin{array}{c}
\\[-6mm]
\scriptstyle L=0 \\[-1mm]
\scriptstyle L+\ell/2 \equiv \pm p \mod{N}
\end{array}}^{\infty} \hspace{-7mm}
a^{L/N} q^{L^2/N} \, \alpha_L,
\end{equation}
where $\ell+\ell'+\sigma N$ must be even.
An alternative form (to that implied by (\ref{e55})) for $c_m^{\ell}$,
known as Hecke's indefinite modular form, reads~\cite{JM,KP}
\begin{eqnarray}
\lefteqn{
c_m^{\ell} = 
\frac{q^{h_m^{\ell}}}{(q)_{\infty}^3} 
\biggl\{ \biggl( \: \sum_{j\geq 0} \sum_{k\geq 0}
- \sum_{j< 0} \sum_{k< 0}\: \biggr)
(-1)^j q^{j(j+\ell+m+1)/2+k((j+k)(N+2)+\ell+1)} \biggr.} \nonumber \\
& &
\hspace{11.4mm} + \biggl. \biggl( \: \sum_{j>0} \sum_{k\geq 0}
- \sum_{j\leq 0} \sum_{k< 0}\: \biggr)
(-1)^j q^{j(j+\ell-m+1)/2+k((j+k)(N+2)+\ell+1)} \biggr\},
\end{eqnarray}
for $\ell=0,\ldots,N-1$, $\ell-m$ even and $|m|\leq \ell$,
with $h_m^{\ell}=\ell(\ell+2)/[4(N+2)]-m^2/(4N)$.
To obtain an expression for $c_m^{\ell}$, when $|m|>\ell$, one
may use the symmetries $c_m^{\ell}=c_{-m}^{\ell}=c_{m+2N}^{\ell}=
c_{N-m}^{N-\ell}$.

\vspace{5mm}
As a final comment and advertisement, we remark that an application
of the higher-level Bailey lemma to obtain $q$-series 
identities for the characters of the level-$N$ A$_1^{(1)}$ 
coset conformal field theories, will be given in ref.~\cite{BMSW}
in collaboration with A.~Berkovich and B.~M.~McCoy.

\section*{Acknowledgements}
We would like to thank G.~E.~Andrews for kind encouragements and 
A.~Berkovich, A.~N.~Kirillov, B.~M.~McCoy, M.~Okado and W.~Orrick 
for valuable discussions and interest in this work. 
We thank D.~L.~O'Brien for helpful suggestions in preparing the manuscript.
This work was (partially)
supported by the Australian Research Council and NSF grant DMR9404747.

\end{document}